\documentclass[10pt]{article}

\usepackage{graphicx}			
\usepackage{amssymb}
\usepackage{amsmath}
\usepackage{verbatim}

\parskip = \baselineskip

\def\vp{\vspace{\baselineskip}}

\title{Matching Rules for the Sphinx Tiling Substitution}
\date{}
\author{Chaim Goodman-Strauss\\ Univ. Arkansas\\ \tt strauss@uark.edu}

\def\tileset{{\mathcal T}}

\begin{document}

\maketitle

This is a copy of notes dated August 14, 2003,  available as~\cite{gs_sphinx_notes}, here 
transliterated into a more traditional format, with some amendments, newly prepared to support the recent {\em Lots of Aperiodic Sets of Tiles}~\cite{gs_LASTs} where the   definitions ``matching rules",  ``enforcing",  a ``substitution tiling", etc. can be found. The construction here is an example of applying the very general techniques of~\cite{gs_mrst}, taking sensible shortcuts as we go. Local symmetries might be found and exploited,  as in every artful example.

The  proof is standard for aperiodic hierarchical tilings~\cite{berger} --- we induct, showing that the tiles clump into larger and larger versions of themselves.

\vp\vp

\centerline{\includegraphics{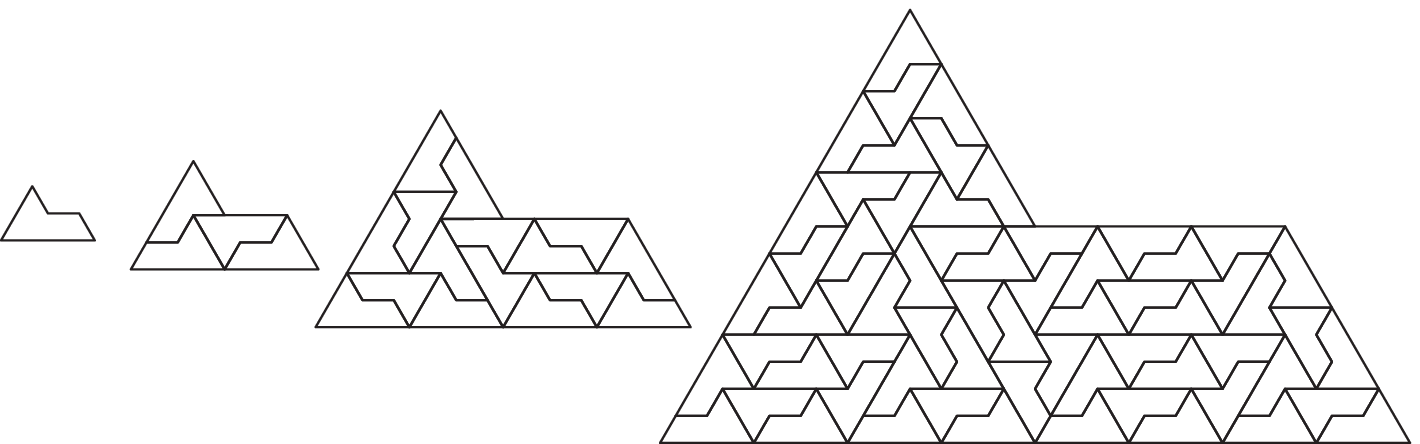}}
\vp
A tiling by sphinx tiles is a sphinx substitution tiling if and only if every finite patch looks like a patch in some inflated sphinx supertile. However the sphinx tiles themselves can tile in lots of different ways. Is there a set of tiles that can {\em only} tile in this manner?

Consider the set $\tileset$ of tiles illustrated on the following pages.  The matching rules are simply that colors match.

{\bf Theorem:} {\em The set $\tileset$ of tiles enforce the sphinx substitution system.}

We encode the structures in~\cite{gs_mrst} as different shades of colors.  Detailed definitions can be found therein.

\newpage
The basic information giving a supertile its immediate identity is its location within its parent (denoted $\mathcal S$ in~\cite{gs_mrst}), shown here in blues:

\centerline{\includegraphics[width=.9\textwidth]{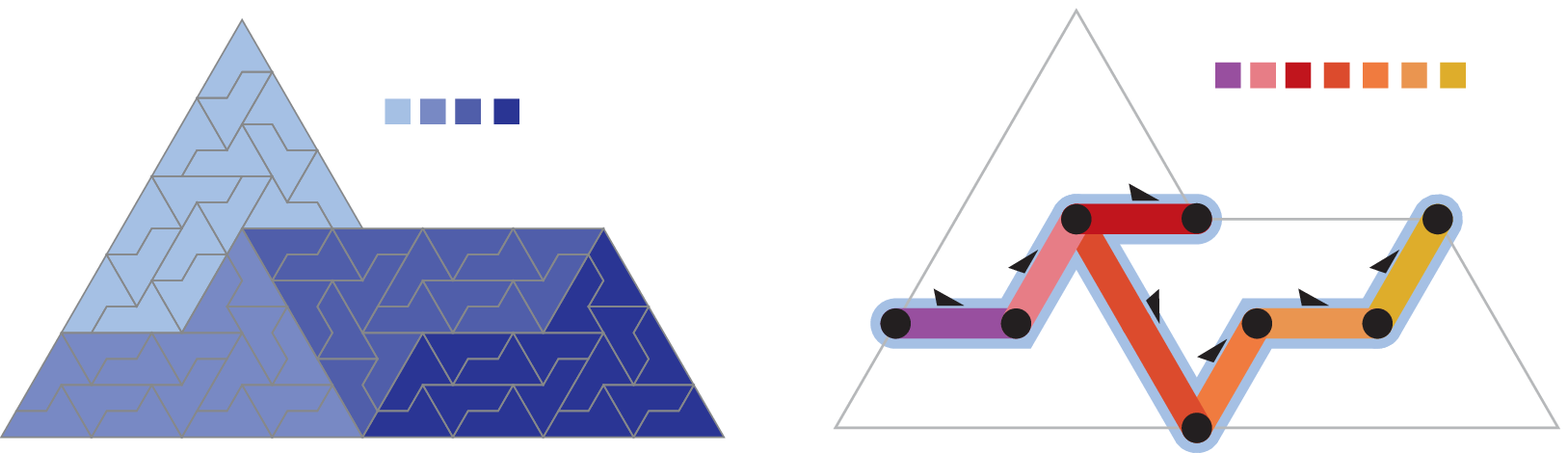}}

Information will be shuttled along a ``skeleton" of edges. Note that every edge in a sphinx substitution tiling lies on the skeleton of a unique sized supertile. (Note in terms of~\cite{gs_mrst}, $\kappa=1$.)

 Each edge then needs to know to what part of a skeleton  it belongs: location, orientation and direction, here encoded in warmer colors and a local framing. We can include a  blue marking, uniform across a skeleton, indication the supertile's position within its parent. 
 
 We store the information particular to a given supertile on the edge tiles on its skeleton. Thus each edge tile has a channel indicating the position of the tile in the skeleton (top), the position of the supertile in its parent (blue) and room for any further information we wish to convey  (in fact, one more channel, for vertex wires, will suffice).
We place an orientation on each  edge tile, partly through its geometry and partly through color. 
This orientation encodes the orientation of each edge of the skeleton. 

We recompose the sphinx tiles, and its supertiles and tilings, into ``edge-", ``vertex-" and ``tile-" tiles, shown at left in the figure below.

\centerline{
\includegraphics{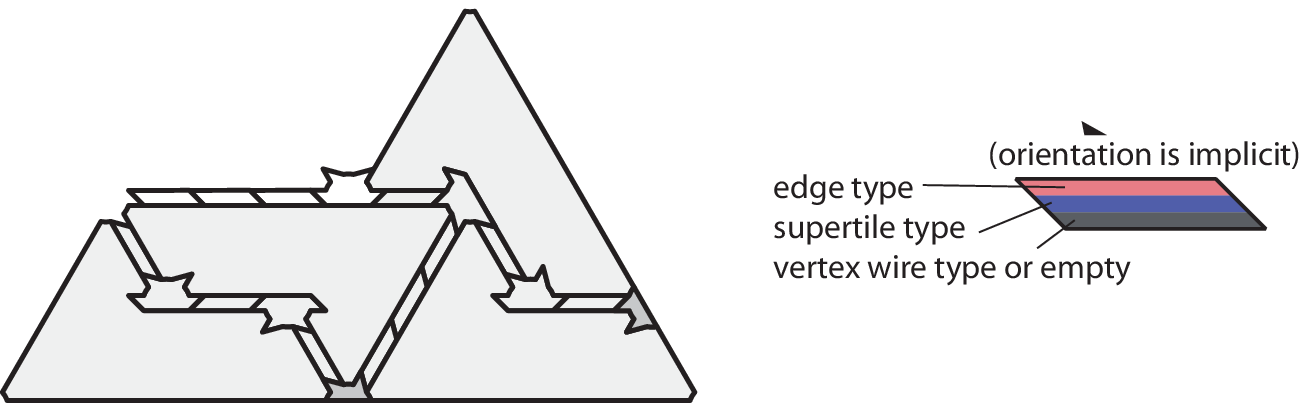}}

The vertex tiles are each peculiar to a specific location in the skeleton, or where a child's skeleton meets its parent skeleton, etc, encoded in its colored markings and geometry. Each sphinx supertile will have a skeleton to convey its identifying information, and the vertices where these markings come together are illustrated below. (Each edge at right has an orientation, as indicated.)

\centerline{
\includegraphics{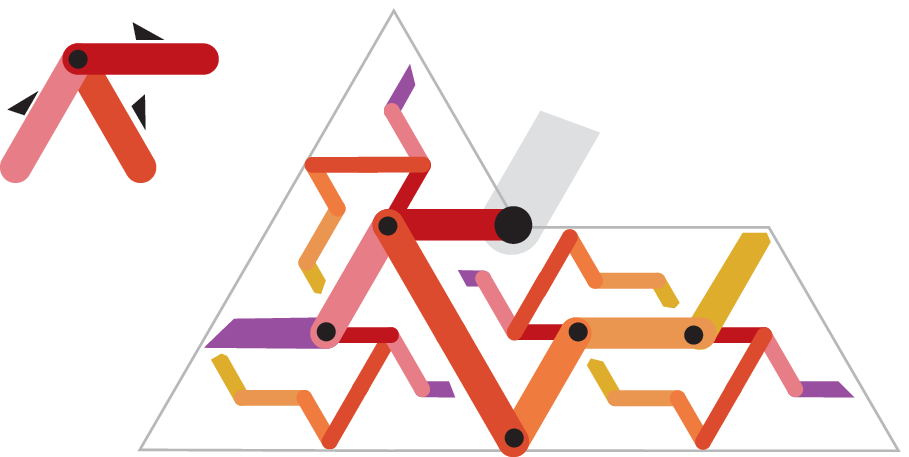}}

Below, we illustrate a typical sphinx vertex tile, which 
mediates between an edge of a child's skeleton (pointing NNW) and  its parent skeleton.  In the 
 first channel we ensure that these skeleton edges   are correctly assembled at this location. In the second channel we specify the identity of the supertile to its parent. (Thus the combinatorics of where children sit within their parent's skeleton is encoded as restrictions on  the second channel of the child with the first channel of the parent.)

\centerline{
\includegraphics[width=\textwidth]{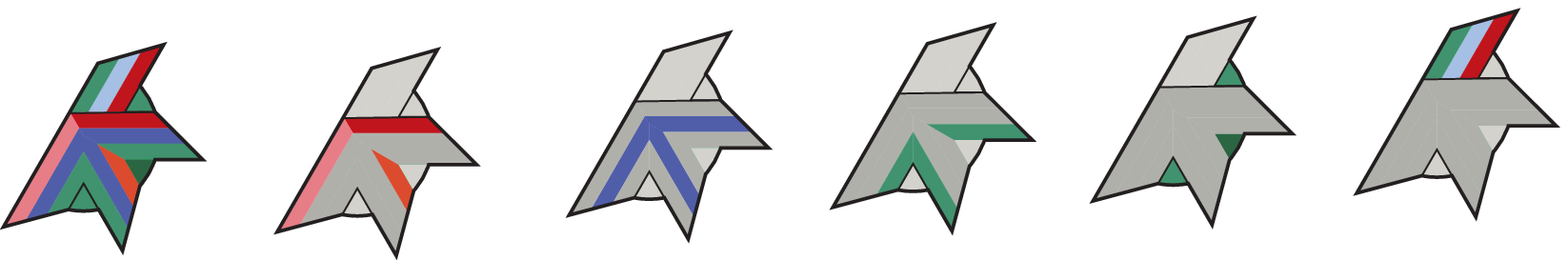}}

However, this is not yet sufficient to {\em ensure} that the tiles enforce the sphinx substitution tiling. Each signal of equivalent edge tiles can propagate indefinitely: we must be sure to place, in just the right spot, a terminate signal (``vertex wire").

In the case of the sphinx substitution, this requires delivering such a signal to any of the vertices shown in green at left below. (In~\cite{gs_mrst}, these are all ``epi-" vertices.)

These vertices themselves form a substitution system: any given green vertex in a parent will be a particular substituted green vertex in the child. In the third channel we place the color of the vertex requiring a terminate signal.

Thus, the vertex tiles further ensure that if a parent is sending a  terminate signal to a particular vertex, so too is the child, and vice versa. 

In this way, we can locate such a vertex within a supertile. As this must be in a specific location on the boundary of its parent, any edge of the parents skeleton meeting this vertex is required to terminate by restrictions on the combinatorics of the vertex tiles. At right we show how this information propogates--- note that for efficiency in defining our edge and vertex tiles, we  put vertex wire information on exactly those edges that might ever be needed to convery it.

\centerline{
\includegraphics[width=\textwidth]{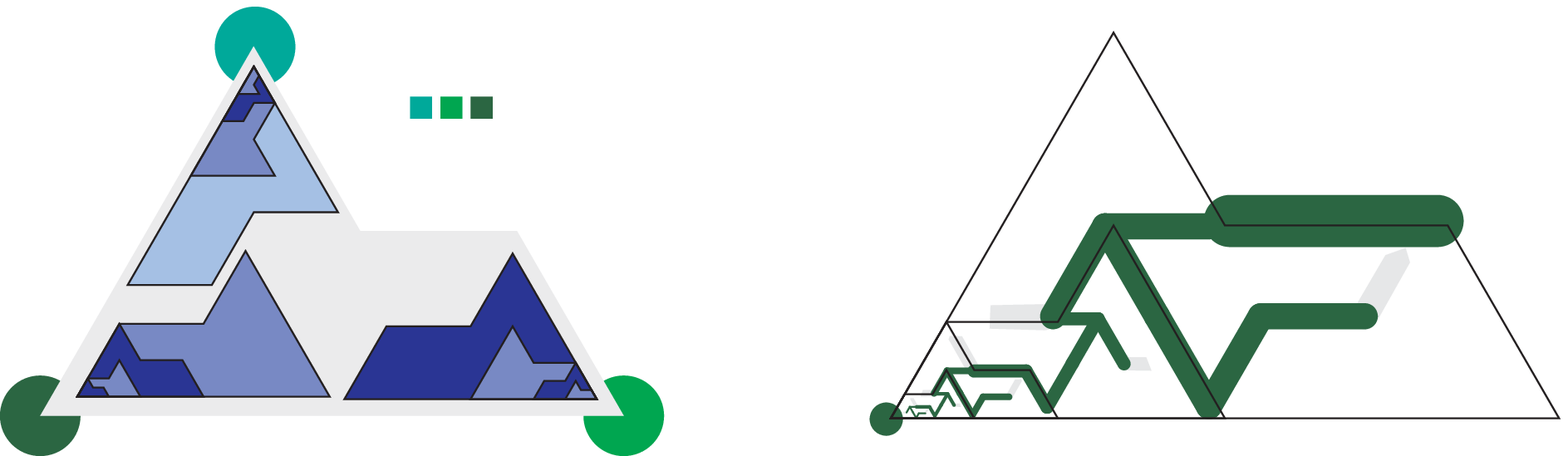}}

A given supertile conveys three pieces of information at its boundary: where it expects to be in its parent, which edges it expects to be there, and passes along any higher level vertex wire information down to its eventual terminal vertex.  

Combinatorially, every supertile formed by these tiles must then be of the form shown below, with information on the boundary, possibly crossed by a  vertex wire, which will originate where the supertile's skeleton meets that of its parent, and terminate at one of two epivertices. One of these is shown:

\centerline{\includegraphics[width=3in]{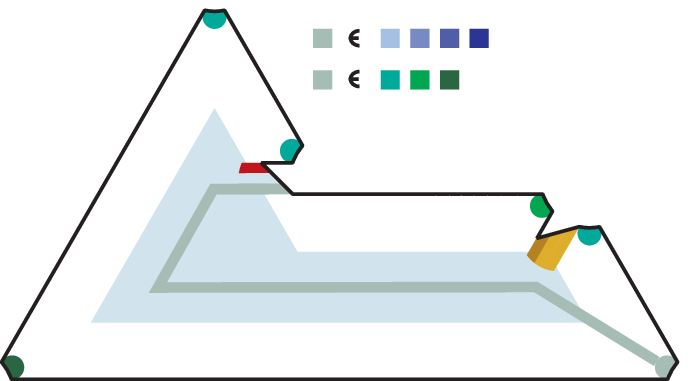}}

\vp\vp 

Taking this as a template for our marked sphinx tiles, we  produce a 1-level supertile. (The presence and color, of vertex wire markings on the two lower marked sphinx tiles depend upon the color of the vertex wire of the parent.) We can see the types of vertex tiles and edge tiles that we need.

\centerline{\includegraphics[width=5in]{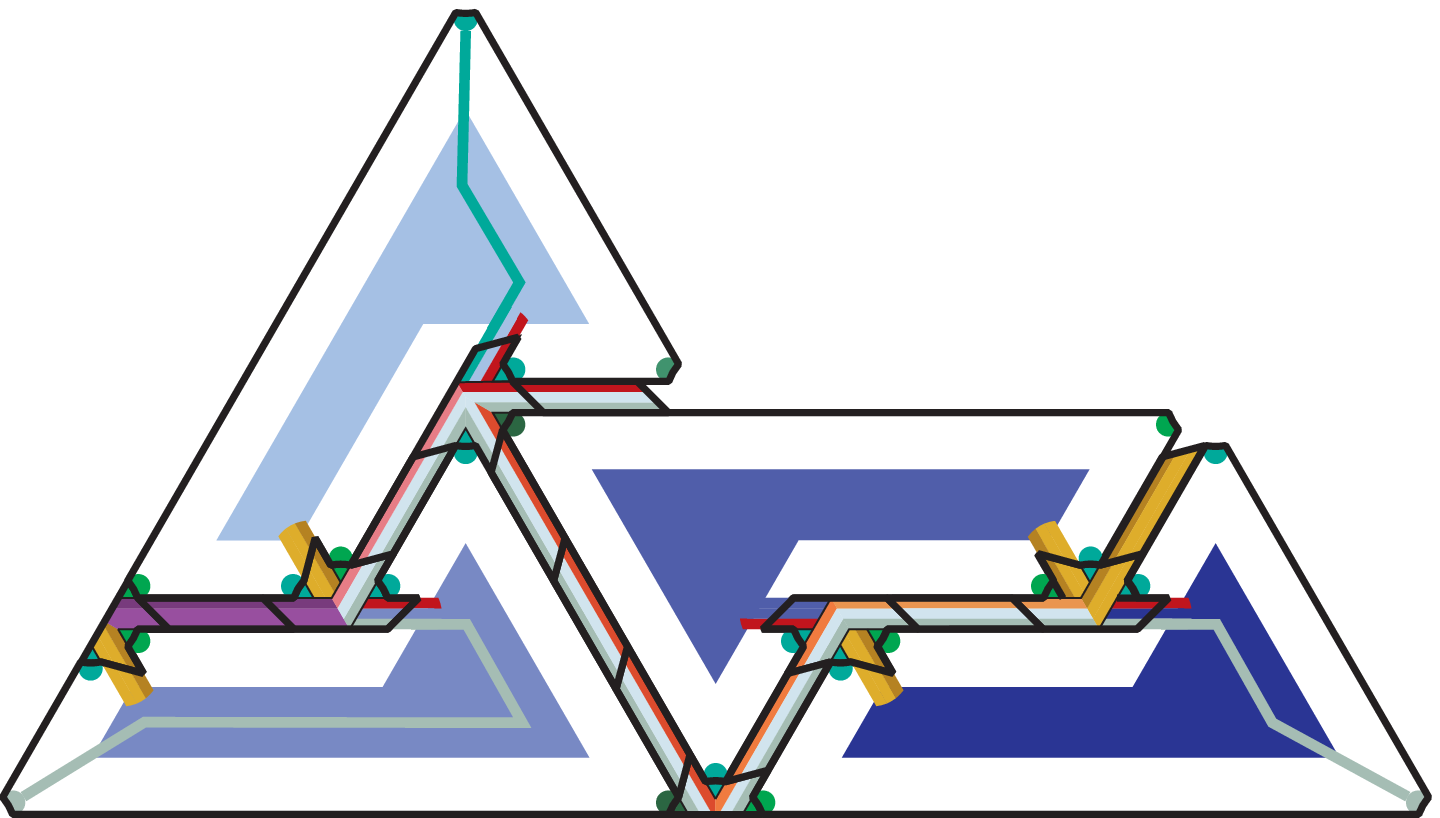}}

And here is a 2nd-level supertile:

\centerline{\includegraphics[width=\textwidth]{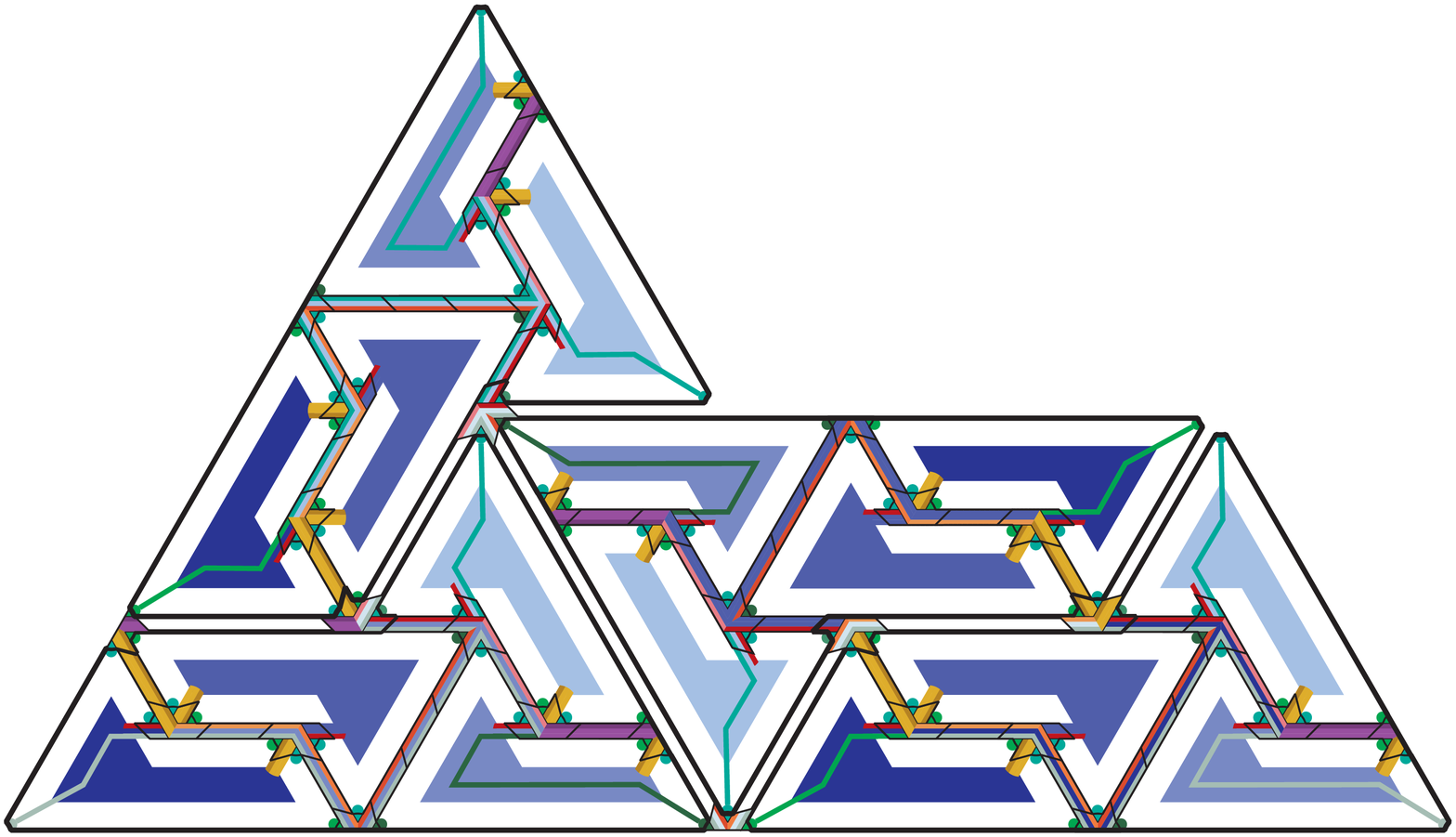}}

We can now easily enumerate the  marked sphinx tiles (and types of supertiles). We only need markings where the supertiles skeleton may meet a higher level supertile:

\centerline{
\includegraphics[scale=.6]{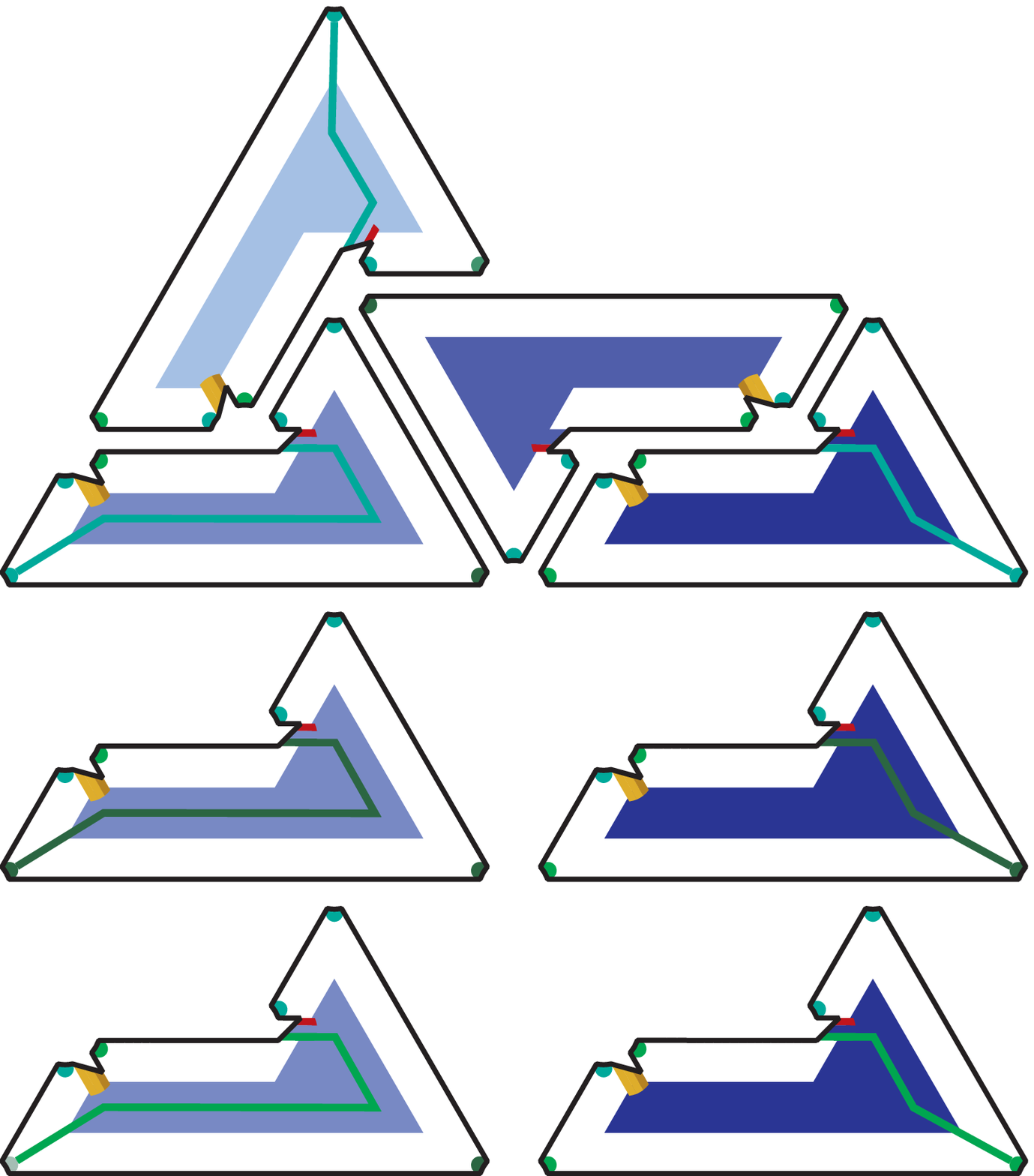}}

\vp

And we enumerate the vertex-, edge-tiles:

\centerline{
\includegraphics[width=\textwidth]{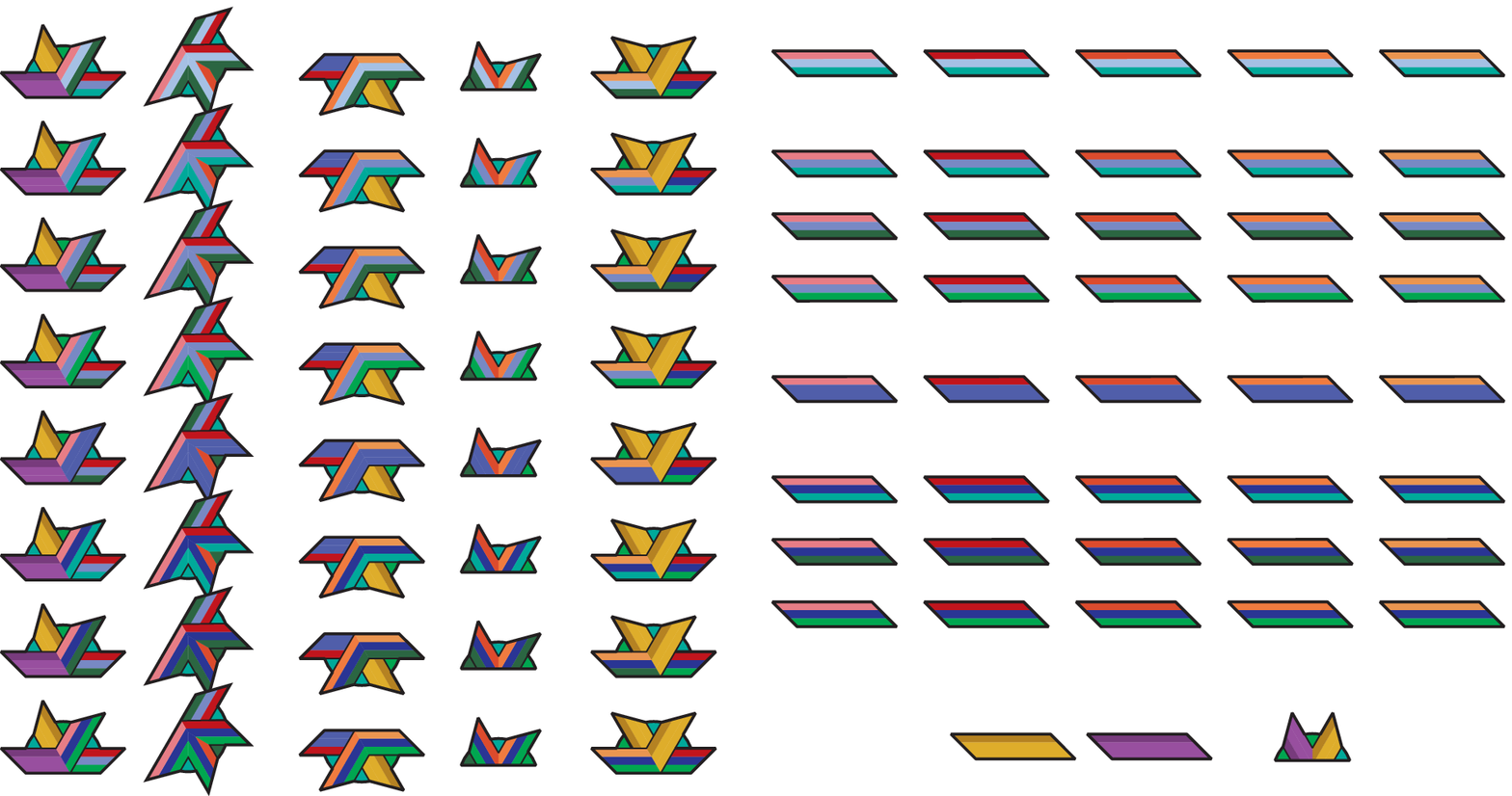}}





\begin{thebibliography}{6}

\bibitem{berger} R. Berger, {\em The undecidability of the domino problem},
Memoirs Am. Math. Soc.  {\bf 66} (1966).

 
\bibitem{gs_mrst} {C. Goodman-Strauss}, {\em Matching rules and substitution
tilings},  Annals of Math. {\bf 147} (1998), 181-223.

%
%

%

%
%
%
%
%
  
  \bibitem{gs_LASTs}{C. Goodman-Strauss} {\em Lots of aperiodic sets of tiles}, arXiv
  
 \bibitem{gs_sphinx_notes}{\em http://comp.uark.edu/~strauss/distribution/tilings/SphinxMatchingRules.pdf}
\end{thebibliography}
\end{document}